\documentclass[11pt]{amsart}

\usepackage{amsmath}
\usepackage{amssymb}
\usepackage{latexsym}
\usepackage{amsfonts}

\begin{document}

\newcommand{\ci}[1]{_{ {}_{\scriptstyle #1}}}

\newcommand{\norm}[1]{\ensuremath{\|#1\|}}
\newcommand{\abs}[1]{\ensuremath{\vert#1\vert}}
\newcommand{\p}{\ensuremath{\partial}}
\newcommand{\pr}{\mathcal{P}}

\newcommand{\wt}{\widetilde}

\newcommand{\pbar}{\ensuremath{\bar{\partial}}}
\newcommand{\db}{\overline\partial}
\newcommand{\D}{\mathbb{D}}
\newcommand{\B}{\mathbb{B}}
\newcommand{\Sp}{\mathbb{S}}
\newcommand{\T}{\mathbb{T}}
\newcommand{\R}{\mathbb{R}}
\newcommand{\C}{\mathbb{C}}
\newcommand{\N}{\mathbb{N}}
\newcommand{\scrH}{\mathcal{H}}
\newcommand{\scrL}{\mathcal{L}}
\newcommand{\td}{\widetilde\Delta}

\newcommand{\la}{\lambda}
\newcommand{\La}{\langle }
\newcommand{\Ra}{\rangle }
\newcommand{\ran}{\operatorname{Ran}}
\newcommand{\tr}{\operatorname{tr}}
\newcommand{\supp}{\operatorname{supp}}
\newcommand{\vf}{\varphi}
\newcommand{\f}[2]{\ensuremath{\frac{#1}{#2}}}

\newcommand{\cJ}{\mathcal{J}}


\newcommand{\entrylabel}[1]{\mbox{#1}\hfill}

\newenvironment{entry}
{\begin{list}{X}%
  {\renewcommand{\makelabel}{\entrylabel}%
      \setlength{\labelwidth}{55pt}%
      \setlength{\leftmargin}{\labelwidth}
      \addtolength{\leftmargin}{\labelsep}%
   }%
}%
{\end{list}}


\numberwithin{equation}{section}

\newtheorem{thm}{Theorem}[section]
\newtheorem{lm}[thm]{Lemma}
\newtheorem{cor}[thm]{Corollary}

\theoremstyle{remark}
\newtheorem{rem}[thm]{Remark}
\newtheorem*{rem*}{Remark}


\title[Carelson Potentials and Embedding Theorems]{Carleson Potentials and the Reproducing Kernel Thesis for Embedding Theorems}

\author[Stefanie Petermichl]{Stefanie Petermichl$^1$}
\address{Stefanie Petermichl, Department of Mathematics\\ University of Texas at Austin\\ 1 University Station C1200\\ Austin, TX USA 78712-0257}
\email{stefanie@math.utexas.edu}
\thanks{1. Research supported in part by a National Science Foundation Grant.}

\author[Sergei Treil]{Sergei Treil$^2$} 
\address{Sergei Treil, Department of Mathematics\\ Brown University\\ 151 Thayer Street Box 1917\\ Providence, RI USA  02912}
\email{treil@math.brown.edu}
\thanks{2. Research supported in part by a National Science Foundation Grant.}

\author[Brett D. Wick]{Brett D. Wick$^3$}
\address{Brett D. Wick, Department of Mathematics\\ Vanderbilt University\\ 1326 Stevenson Center\\ Nashville, TN USA 37240-0001}
\email{brett.d.wick@vanderbilt.edu}
\thanks{3. Research supported in part by a National Science Foundation Grant RTG Grant to Vanderbilt University.}


\begin{abstract}
In this note we present a new proof of the Carleson Embedding Theorem on the unit disc and unit ball in $\C^n$.  The only technical tool used in the proof of this fact is Green's formula.  The starting point is that every Carleson measure gives rise to a bounded subharmonic function.  Using this function we construct a new related Carleson measure that allows for a simple embedding.   In the case of the disc $\D$ this gives the best known constant, with the previous best given by N.~Nikolskii.
\end{abstract}

\maketitle


\setcounter{section}{-1}

\section{Introduction}
The famous Carleson Embedding Theorem for  the unit disc states, in particular, that the embedding of the Hardy space $H^2(\D)$ into a space $L^2(\mu)$ can be checked on reproducing kernels of the Hardy space.  Namely,  it can be stated as follows:

\begin{thm}[Carleson Embedding Theorem]
Let $\mu$ be a non-negative measure in $\D$.  Then the following are equivalent:

\begin{itemize}
\item[(i)] 
The Hardy space $H^2 (\D)$ is embedded in $L^2(\mu)$, i.e.
$$
\int_\D \abs{f(z)}^2 d\mu(z) \le A(\mu)^2 \|f\|_{H^2(\D)}^2 \qquad \forall f\in H^2(\D). 
$$
\item[(ii)] 
\[C(\mu)^{2}:=\sup_{z\in\D}\norm{k_z}^2_{L^{2}(\mu)}=\sup_{z\in\D}\norm{\mathcal{P}_z}_{L^1(\mu)}<\infty,\] 
where $k_z(\xi)=\f{(1-\abs{z}^{2})^{1/2}}{1-\xi\overline{z}}$, is 
the normalized reproducing kernel for the Hardy space $H^{2}(\D)$.
\item[(iii)] 
\[I(\mu):=\sup\left\{\f{1}{r}\mu(\D\cap Q(\xi,r)):r>0,\xi\in\T\right\}<\infty,\] where $Q(\xi,r)$ is a ball in $\C$ with center $\xi$ on $\T$ and radius $r$.
\end{itemize}

Moreover,  the best possible constant $A(\mu)^2$ in (i), the constants $C(\mu)^2$, and $I(\mu)$ are equivalent in the sense of two-sided estimates.
\end{thm}
Property (iii) is typically taken as the definition of a Carleson measure on $\D$. Condition (ii) can be considered as a conformally invariant definition of a Carleson measure. The equivalence  (ii) $\iff$ (iii) above is a simple and standard fact that can be obtained by integrating $|k_z(\cdot)|^2$ using its distribution function. Condition (ii) means that we check the embedding only on the reproducing kernels and not on all $H^2(\D)$ functions. Thus the implication (i) $\implies$ (ii),  as well as the estimate $C(\mu)\le A$ are trivial, so the only non-trivial estimate in this theorem is (ii)$\implies$(i).

The ``reproducing kernel thesis'' is the idea that it is sufficient to check the boundedness of an operator only on reproducing kernels.  The Carleson Embedding Theorem is such an example of this.

In this note we present a new simple proof of the  implication (ii)$\implies$(i) which is quite different from others in the literature.  This proof also gives  the best known estimate of the norm of the embedding operator. Namely, we present a simple proof of the following theorem:

\begin{thm}
\label{CET-Disk}
Suppose
$$
\sup_{\lambda\in\supp \mu}\int_{\D}\abs{k_\lambda(z)}^2d\mu(z)=: A <\infty. 
$$
Then
$$
\int_{\D}\abs{f(z)}^2d\mu(z)\le  2e A\|f\|^2_{H^2(\D)},\qquad \forall f\in H^2. 
$$
\end{thm}

Recall that the Hardy space can be defined as the closure of analytic polynomials in $L^2(\T, m)$, where $m$ is normalized ($m(\T)=1$) Lebesgue measure on $\T=\p\D$ with the norm inherited from $L^2(\T, m)$. The elements of $H^2(\D)$ admit a natural analytic continuation inside the unit disc $\D$ (see \cite{Garnett}), so the integral $\int_{\D}\abs{f(z)}^2d\mu(z)$ in the above theorem is defined. 

Note that the theorem says that
 it is sufficient to check the embedding not on all reproducing kernels $k_\la$, but only on $k_\la$, $\la\in \supp\mu$. This fact was  known before, cf. \cite[p.~151]{Nik2}, but the constant $2e$ is the best known to date.  In \cite[p.~151]{Nik2} the constant $32$ in the reproducing kernel thesis was obtained, and later in \cite[p.~105]{Nik1} the constant was improved to $16$.

The proof we are going to present is a simple ``conformally invariant'' proof with the main tool used being Green's formula. This proof generalizes easily to the unit ball in $\C^n$. 

Recall, that the Hardy space $H^2(\B_n)$ on the unit ball in $\C^n$ can be defined as the closure of polynomials in $L^2(\Sp, \sigma)$, where $\sigma$ is  the Lebesgue measure  on the boundary $\Sp_n=\p\B_n$ (see \cite{Krantz}, or \cite{Rudin} for other equivalent definitions, as well as for more information about this space).

We obtain the following ``reproducing kernel thesis'' for $H^2(\B_n)$. 
\begin{thm}
\label{CET-Ball}
Suppose
$$
\sup_{\lambda\in\supp \mu}\int_{\B_n}\abs{k_\lambda(z)}^2d\mu(z)= C , 
$$
where $k_\la $, $\la \in \B_n$, is the normalized reproducing kernel of $H^2(\B_n)$. 
Then
$$
\int_{\B_n}\abs{f(z)}^2d\mu(z)\le  e\frac{(2n)!}{(n!)^2} C\|f\|^2_{H^2(\B_n)}, \qquad \forall f\in H^2(\B_n).
$$
\end{thm}

\begin{rem}
The statement of the theorem does not depend on the choice of normalization of the measure $\sigma$ because if one replaces $\sigma $ by $c\sigma$ one would need to multiply the reproducing kernel by $c^{-1/2}$. Usually normalization is chosen by assuming that $\sigma(\Sp)=1$ and in this case the reproducing kernel $k_\lambda$ is given by (see \cite{Rudin})
$$
k_\lambda(z)=\f{(1-\abs{\lambda}^{2})^{n/2}}{(1-\La z,\la \Ra)^n}, 
$$
with $\La\cdot,\cdot\Ra$ denoting the standard Hermitian inner product in $\C^n$.
\end{rem}

\begin{rem}
The theorem is well known, and is usually proved by real variable methods. A new part here would be the estimate and the fact that it is  sufficient to check the embedding only on $k_\la$, $\la\in \supp \mu$.  We do not see how to immediately get the latter from known results, short of repeating the proof of the Carleson Embedding Theorem given in \cite{Nik1} in the context of the unit ball. 
\end{rem}
The authors would like to thank Alexander Volberg  and Dechao Zheng for  useful discussions.

Throughout the paper the notation $:=$ means equal by definition, and $A\lesssim B$ means there exists an absolute positive constant $C$ such that $A\leq CB$.  The expression $A\approx B$ means $A\lesssim B$ and $B\lesssim A$ both hold.


\section{The Embedding Theorem for the Unit Disc $\D$}

\subsection{Uchiyama's Lemma}
We need the following Lemma, a version of which was probably first proved by Uchiyama.
\begin{lm}
\label{Uchi}
Suppose that $\varphi\le 0$ is a subharmonic function.  Then 
$$d\nu(z):=\frac{1}{2\pi} e^{\varphi}\Delta\varphi(z) \log\f{1}{\abs{z}}dA(z)$$ is a Carleson measure and the embedding $H^2(\D)\subset L^2(\nu)$ is a contraction.  More precisely, for any $f\in H^2(\D)$ we have
$$
\int_{\D}\abs{f(z)}^2d\nu(z)\leq \norm{f}_{H^2(\D)}^2.
$$
\end{lm}
\begin{proof}
To prove this lemma we will simply use Green's Formula applied to a particular function.  First, recall that Green's Formula for a function $u$ says
$$
\f{1}{2\pi}\int_{\D}\Delta u(z)\log\f{1}{\abs{z}}dA(z)=\int_{\T}u(\xi)dm(\xi)-u(0),
$$
where $m$ is the normalized ($m(\T)=1$) Lebesgue measure on the unit circle $\T=\p \D$. 

We now let $u=e^\varphi\abs{f}^2$.  Let us  compute the Laplacian of this function.  Recalling  the definition of $\p$- and $\pbar$-derivatives, 
$$
\p f = \frac12 \left( \frac{\p f}{\p x} - i \frac{\p f}{\p y} \right), \qquad 
\pbar f = \frac12 \left( \frac{\p f}{\p x} + i \frac{\p f}{\p y} \right)
$$
and the fact that $\Delta= 4\p \pbar$ we get
\begin{equation}
\label{1.1}
\Delta (e^\varphi\abs{f}^2)=e^\varphi(\Delta\varphi \abs{f}^2+ 4\abs{(\p\varphi f+\p f)}^2)\geq e^\varphi\Delta\varphi\abs{f}^2.
\end{equation}
Applying the information of $\varphi$ we have the right hand side of Green's Formula giving
$$
\int_{\T}e^{\varphi(\xi)}\abs{f(\xi)}^2dm(\xi) -e^{\varphi(0)}\abs{f(0)}^2\leq \int_{\T}\abs{f(\xi)}^2dm(\xi).
$$
On the other hand, we have
$$
\f{1}{2\pi}\int_{\D}\Delta(e^{\varphi(z)}\abs{f(z)}^2)\log\f{1}{\abs{z}}dA(z)\geq \f{1}{2\pi}\int_{\D}e^{\varphi(z)}\Delta\varphi(z) \abs{f(z)}^2\log\f{1}{\abs{z}}dA(z).
$$
Combining things we find that 
$$
\f{1}{2\pi}\int_{\D}e^{\varphi(z)}\Delta\varphi(z) \abs{f(z)}^2\log\f{1}{\abs{z}}dA(z)\leq \int_{\T}\abs{f(\xi)}^2dm(\xi)
$$
which gives the Lemma and shows that $e^{\varphi(z)}\Delta\varphi(z)\log\f{1}{\abs{z}}dA(z)$ is a Carleson measure on $\D$.
\end{proof}

\begin{cor}
\label{disketophi}
If $\vf$ is bounded (and we still assume that $\varphi\leq 0$) then $d\nu(z):=\frac{1}{2\pi} \Delta\varphi(z) \log\f{1}{\abs{z}}dA(z)$ is a Carleson measure and for any $f\in H^2(\D)$ we have
$$
\frac{1}{2\pi}\int_{\D}\abs{f(z)}^2\Delta\varphi(z)\log\f{1}{\abs{z}}dA(z)\leq e \|\vf\|_\infty \norm{f}_{H^2(\D)}^2.
$$
\end{cor}
\begin{proof}
Since $\vf \ge -r:= - \|\vf\|_\infty$, Uchiyama's Lemma (Lemma \ref{Uchi}) implies that 
$$
e^{-r} \int_\D \abs{f(z)}^2 d\nu \le \|f\|_{H^2(\D)}^2. 
$$
Replacing $\vf$ by $t\vf$, $t>0$ we get
$$
te^{-tr} \int_\D \abs{f(z)}^2 d\nu \le \|f\|_{H^2(\D)}^2. 
$$
The function $te^{-tr}$ attains its maximum at  $t= 1/r = 1/\|\vf\|_\infty$. Plugging in this value of $t$ we get  the desired estimate.
\end{proof}

\subsection{Carleson Potentials and the Proof of Theorem \ref{CET-Disk}}
Suppose the measure $\mu$ satisfies the assumption of Theorem \ref{CET-Disk}. By homogeneity we can assume without loss of generality that the constant $C$ is $1$.

Define the Carleson Potential
$$
\varphi(z):=-\int_{\D}\abs{k_z(\lambda)}^2d\mu(\lambda)=-\int_{\D}\mathcal{P}_z(\lambda)d\mu(\lambda),
$$
where $k_z$ is the (normalized) reproducing kernel and $\mathcal{P}_z(\lambda)=|k_z(\la)|^2$ is the Poisson kernel at $z$.
Then  $-1\le \varphi(z)\le 0$ for $z\in \supp \mu$.   

We next compute the Laplacian of the function $\varphi(z)$.  Using the fact that for an analytic function $f$ we have $\Delta |f|^2 = \p \pbar |f|^2 = 4 |f'|^2$ we get 
$$
\Delta_z\mathcal{P}_z(\lambda)=4\f{\abs{\lambda}^2-1}{\abs{1-\overline{\lambda}z}^4},
$$
(here $\Delta_z$ stands for the Laplacian in the variable $z$). 
This clearly implies that $\vf$ is subharmonic and that 
$$
\Delta \vf(z) = 4 \int_\D \frac{1-|\la|^2}{|1-\overline \la z|^4}d\mu(\la).
$$

Applying Uchiyama's Lemma (Lemma \ref{Uchi}) we get  
$$
\int_{\D}\abs{f(z)}^2d\nu(z)\leq \norm{f}_{H^2(\D)}^2,
$$
with $d\nu(z):=e^{\varphi(z)}\Delta\varphi(z) \log\f{1}{\abs{z}}dA(z)$.

We will now prove the estimate
\begin{equation}
\label{1.2}
\int_{\D}\abs{f(\la)}^2d\mu(\la)\le 2 e \int_{\D}\abs{f(z)}^2d\nu(z)
\end{equation}
which will immediately imply the theorem.

First note that
$$
\int_{\D}\abs{f(z)}^2d\nu(z)=\f{4}{2\pi}\int_{\D}\int_{\D}\abs{f(z)}^2e^{\varphi(z)}\f{1-\abs{\lambda}^2}{\abs{1-\overline{\lambda}z}^4}\log\f{1}{\abs{z}}dA(z)d\mu(\lambda).
$$
Using  the estimate $\f{1}{2}(1-\abs{z}^2)\leq\log\f{1}{\abs{z}}$ we have 

$$
\int_{\D}\abs{f(z)}^2d\nu(z)\geq\f{1}{\pi}\int_{\D}\int_{\D}\abs{f(z)}^2e^{\varphi(z)}\f{(1-\abs{\lambda}^2)(1-\abs{z}^2)}{\abs{1-\overline{\lambda}z}^4}dA(z)d\mu(\lambda).
$$
\begin{rem} If we did not care about the constant then Theorem \ref{CET-Disk} would be proved.  Here is why.  In the disc centered at $\lambda$ of radius $\f{\delta}{10}>0$ where $\delta=\textnormal{dist}(\lambda, \T)$, call it $D(\lambda,\delta)$,  we have that
$$
\f{1-\abs{\lambda}^2}{\abs{1-\overline{\lambda}z}^4}(1-\abs{z}^2)\approx\f{1}{\delta^2}.
$$
Using the subharmonicity of $e^\vf |f|^2$ (see \eqref{1.1}) and the trivial fact that the volume of $D(\lambda,\delta)\approx\delta^2$ we get 
$$
e^{\varphi(\lambda)}\abs{f(\lambda)}^2\lesssim\int_{D(\lambda,\delta)}e^{\varphi(z)}\abs{f(z)}^2  
\f{(1-\abs{\lambda}^2)(1-\abs{z}^2)}{\abs{1-\overline{\lambda}z}^4}  dA(z) 
$$
Increasing the domain of integration to the whole disc $\D$ clearly does not spoil the inequality, and integrating both sides with respect to $d\mu(\la)$ we obtain   
$$
\int_{\D}e^{\varphi(\lambda)}\abs{f(\lambda)}^2d\mu(\lambda)\lesssim \int_{\D}\abs{f(z)}^2d\nu(z) \lesssim \|f\|^2_{H^2(\D)}
$$
which proves the theorem (without constants).
\end{rem}

However, since we are after the  constants, here is how to obtain the sharper estimate.  We focus on the inner integral and will prove the inequality  
\begin{equation}
\label{1.3}
\f{1}{\pi}\int_{\D}\abs{f(z)}^2e^{\varphi(z)}\f{(1-\abs{\lambda}^2)(1-\abs{z}^2)}{\abs{1-\overline{\lambda}z}^4}dA(z)\geq \f{1}{2}e^{\varphi(\lambda)}\abs{f(\lambda)}^2\quad\forall\lambda\in\textnormal{supp}\ \mu, 
\end{equation}
which after integration  with respect to $d\mu(\la)$ gives \eqref{1.2}.  

Let $w=b_\lambda(z):=\f{\lambda-z}{1-\overline{\lambda}z}$ denote a conformal change of variables (note that $z=b_\la(w)$).  A simple computation shows that
$$
dA(w)=\left(\f{1-\abs{\lambda}^2}{\abs{1-\overline{\lambda}z}^2}\right)^2dA(z).
$$
If we let $\tilde{g}(w):=g\circ b_\lambda(w)$ then the above integral can be recognized as
$$
\f{1}{\pi}\int_{\D}\abs{f(z)}^2e^{\varphi(z)}\f{(1-\abs{\lambda}^2)(1-\abs{z}^2)}{\abs{1-\overline{\lambda}z}^4}dA(z)=\f{1}{\pi}\int_{\D}e^{\tilde{\varphi}(w)}\abs{\tilde{f}(w)}^2\f{1-\abs{w}^2}{\abs{1-\overline{\lambda}w}^2}dA(w).
$$
In this reduction we have used the algebraic identity that for $b_\lambda$ defined above,
$$
1-\abs{z}^2=\f{(1-\abs{\lambda}^2)(1-\abs{w}^2)}{\abs{1-\overline{\lambda}w}^2}.
$$

Continuing the estimate we have 
$$
\f{1}{\pi}\int_{\D}e^{\tilde{\varphi}(w)}\abs{\tilde{f}(w)}^2\f{1-\abs{w}^2}{\abs{1-\overline{\lambda}w}^2}dA(w)=\f{1}{\pi}\int_{\D}e^{\tilde{\varphi}(w)}\left\vert\f{\tilde{f}(w)}{1-\overline{\lambda}w}\right\vert^2(1-\abs{w}^2)dA(w).
$$
The function $\f{\tilde{f}(w)}{1-\overline{\lambda}w}$ is analytic and $\tilde \varphi$ is subharmonic, so (see \eqref{1.1}) the function $u(w) = e^{\tilde{\varphi}(w)}\left| \f{\tilde{f}(w)}{1-\overline{\lambda}w}\right|^2$ is subharmonic.   

Integrating in polar coordinates and using the mean value property for subharmonic functions we get 
$$
\int_\D u(w) (1-|w|^2 ) d A(w) = \int_0^1 (1-r^2)r \int_0^{2\pi} u(r\theta) d\theta dr 
\ge 2\pi u(0) \int_0^1 (1-r^2)r dr = \frac{\pi}{2} u(0).
$$

Gathering all together we find 
$$
\f{1}{\pi}\int_{\D}e^{\tilde{\varphi}(w)}\abs{\tilde{f}(w)}^2\f{1-\abs{w}^2}{\abs{1-\overline{\lambda}w}^2}dA(w)\geq \frac12 e^{\tilde{\varphi}(0)}\abs{\tilde{f}(0)}^2= \f{1}{2}e^{\vf(\la)}\abs{f(\lambda)}^2.
$$
which is equivalent to \eqref{1.3}. 

This finally shows that for a Carleson measure $\mu$ on $\D$ we have
$$
\int_{\D}\abs{f(z)}^2d\mu(z)\leq 2e\norm{\varphi}_{\infty}\norm{f}_{H^2(\D)}^2=2e\norm{\mu}_{\mathcal{C}}\norm{f}_{H^2(\D)}^2
$$
proving Theorem \ref{CET-Disk} for the disc $\D$. \hfill\qed 

 We should also say that the constant $2e$ is the best known constant obtained for the norm of the embedding operator.  N.~Nikolskii has a different proof of the Carleson Embedding Theorem in which the constant obtained is $32$.  See either \cite{Nik1} or \cite{Nik2} for the proof.  We further conjecture that the constant $e$ is sharp in Uchiyama's Lemma (Lemma \ref{Uchi}) and the constant $2e$ is sharp in Theorem \ref{CET-Disk}.

We will use the proof in the disc as motivation for the appropriate proof on the unit ball in $\C^n$.


\section{The Embedding Theorem for the Unit Ball $\B_n$}

The proof of Theorem \ref{CET-Ball} is very similar to the one given for the case of the disc $\D$.  The essential difference is that  one must use the invariant Laplacian for the unit ball instead of the usual Laplacian.  This reflects  the complex structure of the unit ball $\B_n$.  In particular, the embedding theorem is usually stated in terms of  ``Carleson cubes''  defined via the non-isotropic metric, as opposed to the standard euclidean one.  The other motivation for the use of invariant Laplacian follows from the fact that $|k_\la|^2$ is the invariant Poisson kernel. 

Recall that the invariant Laplacian is defined by the following formula
$$
\tilde{\Delta}:=4\sum_{i,j}g^{ij}\pbar_i\partial_j
$$
with $g^{ij}=\f{1-\abs{z}^2}{n+1}(\delta_{ij}-\overline{z}_iz_j)$ the components of the inverse of the Bergman metric on $\B_n$, and  
$$
\p_j f = \frac12 \left( \frac{\p f}{\p x_j} - i \frac{\p f}{\p y_j} \right), \qquad 
\pbar_j f = \frac12 \left( \frac{\p f}{\p x_j} + i \frac{\p f}{\p y_j} \right)\quad\forall j=1,\ldots,n.
$$

We first need to translate Uchiyama's Lemma to the ball.  

\subsection{Uchiyama's Lemma for the Unit Ball}
We need the following variant of Lemma \ref{Uchi}.  The appropriate analog of Uchiyama's Lemma on the ball requires a few minor modifications to deal with the additional number of variables and the complex structure.  We use the Green's Function (with the pole at 0) for the invariant Laplacian, which is given by
$$
G(\lambda)=\f{n+1}{2n}\int_{\abs{\lambda}}^1(1-t^2)^{n-1}t^{-2n+1}dt,
$$
and in the case $n=1$ reduces to the usual logarithm.  This function will play the same role that the logarithm plays in the disc.  See \cite{Stoll} for the derivation of Green's function $G(\lambda)$ for the invariant Laplacian.  We also need to use the volume form, or the invariant measure on the unit ball.  It is given by
$$
dg(\lambda):=\f{dV(\lambda)}{(1-\abs{\lambda}^2)^{n+1}}
$$ 
with $dV$ the standard (non-normalized) volume form for the unit ball $\B_n$.

\begin{lm}
\label{Uchi-Ball}
Suppose that $\varphi$ is a non-positive invariant subharmonic function.  Then 
$$d\nu(z):=\frac{n!}{\pi^n} e^{\varphi(z)}\widetilde{\Delta}\varphi(z)G(z)dg(z)$$ is a Carleson measure and the embedding $H^2(\B_n)\subset L^2(\nu)$ is contractive.  More precisely, for any $f\in H^2(\B_n)$ we have
$$
\int_{\B_n}\abs{f(z)}^2d\nu(z)\leq \norm{f}_{H^2(\B_n)}^2.
$$
\end{lm}

\begin{cor}
\label{balletophi}
If $\vf$ is bounded (and we still assume that $\varphi\leq 0$) then $d\nu(z):=\frac{n!}{\pi^n} \widetilde{\Delta}\varphi(z) G(z)dg(z)$ is a Carleson measure and for any $f\in H^2(\B_n)$ we have
$$
\frac{n!}{\pi^n}\int_{\B_n}\abs{f(z)}^2\widetilde{\Delta}\varphi(z)G(z)dg(z)\leq e \|\vf\|_\infty \norm{f}_{H^2(\B_n)}^2.
$$
\end{cor}

\begin{proof}[Proof of Lemma \ref{Uchi-Ball}]
 We begin by showing that
\begin{equation}
\label{2.1}
\widetilde{\Delta}\left(e^{\varphi}\abs{f}^2\right) \ge (\widetilde{\Delta} \varphi )|f|^2 e^{\varphi}.
\end{equation}
Indeed, using the chain rule and that $f$ is a holomorphic function we arrive at
\begin{eqnarray*}
\lefteqn{\widetilde{\Delta}\left(e^{\varphi}\abs{f}^2\right)}\\ 
& = & 4\sum_{i,j}g^{ij}\left[\pbar_i\varphi\p_j\varphi\abs{f}^2+\p_j\varphi f\overline{\p_i f}+\pbar_i\varphi\p_j f\overline{f}+\p_j f\overline{\p_i f}+\pbar_i\p_j\varphi\abs{f}^2\right]e^\varphi\\
& = & (\widetilde{\Delta}\varphi )\abs{f}^2 e^\varphi +4 \sum_{i,j}g^{ij}\La\p_i f+\p_i\varphi f,\p_j f+\p_j\varphi f\Ra e^\varphi\\
& = & \widetilde{\Delta}\varphi\abs{f}^2 e^\varphi + 4e^\varphi\norm{\p f+\p\varphi f}^2_{Berg}\\
&\ge &(\widetilde{\Delta} \varphi )|f|^2 e^{\varphi}.
\end{eqnarray*}

The rest of the proof is the standard application of the Green's formula. Green's formula in the Bergman metric is given by (see \cite{Chang} or \cite{Zheng})
$$
\f{n!}{\pi^n}\int_{\B}\widetilde{\Delta}u(z) G(z) dg(z) = \int_{\Sp} u(\xi) d\sigma(\xi) -u(0).
$$
where $d\sigma$ is the normalized Lebesgue measure for $\Sp$, i.e., $\sigma(\Sp)=1$.  The formula in \cite{Zheng} is derived, but the exact constants weren't computed, however since we are after sharp constant this more precise formula is important.  The precise constants can be derived from \cite{Zheng} by testing Green's formula on the radial function $f(z)=1-\abs{z}^2$ and then performing straight forward, though tedious, computations.

So using \eqref{2.1} and applying Green's  formula with $u=e^\varphi\abs{f}^2$ we continue our estimate
\begin{eqnarray*}
\int_{\B_n} |f(z)|^2 d\nu(z) \le\f{n!}{\pi^n}\lefteqn{\int_{\B} \widetilde{\Delta}_{z} (e^{\varphi(z)}|f(z)|^2) G(z) dg(z)}\\
& = & \int_{r\Sp} e^{\varphi(\xi)}\abs{f(\xi)}^2 d\sigma(\xi) -C(n)e^{\varphi(0)}\abs{f}^2(0)\\
& \leq & \int_{r\Sp}\abs{f(\xi)}^2d\sigma(\xi).
\end{eqnarray*}
\end{proof}

The proof  Corollary \ref{balletophi} is exactly the same as in the case of the disc, and we leave it as an exercise for  the reader.

\subsection{Carleson Potentials and the Proof of Theorem \ref{CET-Ball}}
We again suppose the measure $\mu$ satisfies the assumption of Theorem \ref{CET-Ball}.  By homogeneity we can assume without loss of generality that the constant $C$ is $1$. 

Define the Carleson Potential
$$
\varphi(z):=-\int_{\B_n}\abs{k_z(\lambda)}^2d\mu(\lambda)=-\int_{\B_n}\mathcal{P}_z(\lambda)d\mu(\lambda),
$$
where $k_z$ is the (normalized) reproducing kernel and $\mathcal{P}_z(\lambda)=|k_z(\la)|^2$ is the Poisson kernel at $z$ for the unit ball $\B_n$, i.e. $\mathcal{P}_z(\lambda)=\f{(1-\abs{z}^2)^n}{\abs{1-\La\lambda, z\Ra}^{2n}}$.  Then  $-1\le \varphi(z)\le 0$ for $z\in \supp \mu$.  The following lemma will be important in computing the invariant Laplacian of the Carleson potential $\vf$.

\begin{lm}
\label{Lap}
Let $\mathcal{P}_z(\lambda)$ denote the Poisson-Szeg\"o kernel.  Then
$$
\widetilde{\Delta}_{z}\mathcal{P}_{z}(\lambda)=-\f{4n^2}{n+1}(1-\abs{z}^2)\mathcal{P}_{z}(\lambda)\mathcal{P}_\lambda(z)^{1/n}.
$$
\end{lm}
It is clear that this Lemma implies that $\varphi$ is invariant subharmonic because upon passing the invariant Laplacian inside the integral we are left with 
$$\widetilde{\Delta}_{z}\varphi(z)=\f{4n^2}{n+1}(1-\abs{z}^2)\int_{\B}\mathcal{P}_{z}(\lambda)\mathcal{P}_\lambda(z)^{1/n}d\mu(\lambda)\geq 0,
$$ 
which is the characterization of (smooth) invariant subharmonic functions. 

\begin{proof}
The proof of this lemma is a straightforward, though tedious, computation.  A simple computation shows
$$
\p_j\pr_{z}(\lambda)=n\left[\f{\overline{\lambda}_j}{1-\La z,\lambda\Ra}-\f{\overline{z}_j}{1-\abs{z}^2}\right]\pr_z(\lambda).
$$
Using that $\pr_z(\lambda)$ is real valued and $\overline{\p_j H}=\pbar_j \overline{H}$ for any function $H$ we have
$$
\pbar_j\pr_{z}(\lambda)=n\left[\f{\lambda_j}{1-\La \lambda,z\Ra}-\f{z_j}{1-\abs{z}^2}\right]\pr_z(\lambda).
$$
Combining this we find that
\begin{eqnarray*}
\lefteqn{\pbar_i\p_j\pr_{z}(\lambda)}\\ 
& = & 
n\pr_{z}(\lambda)\left[n\left(\f{\overline{\lambda}_j}{1-\La z,\lambda\Ra}-\f{\overline{z}_j}{1-\abs{z}^2}\right)\left(\f{\lambda_i}{1-\La \lambda,z\Ra}-\f{z_i}{1-\abs{z}^2}\right)-\left(\f{\delta_{ij}}{1-\abs{z}^2}+\f{z_i\overline{z}_j}{(1-\abs{z}^2)^2}\right)\right].
\end{eqnarray*}
Now, by definition
$$
\widetilde{\Delta}_{z}\pr_{z}(\lambda)=\frac{4}{n+1}(1-\abs{z}^2)\sum_{i,j}(\delta_{ij}-\overline{z}_iz_j)\pbar_i\p_j\pr_{z}(\lambda)
$$
where $\delta_{ij}$ is the Kronecker delta function.  If one is patient enough, then computation yields
$$
\widetilde{\Delta}_{z}\mathcal{P}_{z}(\lambda)=-\f{4n^2}{n+1}(1-\abs{z}^2)\mathcal{P}_{z}(\lambda)\mathcal{P}_\lambda(z)^{1/n}.
$$
\end{proof}

This computation can also been seen by noting that for a K\"ahler manifold, we have for an analytic function $f$ that $\widetilde{\Delta} |f|^2 =  4 |\widetilde{\nabla}f|^2$, where $\widetilde{\nabla}$ denotes the invariant gradient associated to the Bergman metric.  See  \cite{Stein} or \cite{Stoll}.

Applying Uchiyama's Lemma (Lemma \ref{Uchi-Ball}) we get  
$$
\int_{\B_n}\abs{f(z)}^2d\nu(z)\leq \norm{f}_{H^2(\B_n)}^2,
$$
with $d\nu(z):=\f{n!}{\pi^n}e^{\varphi(z)}\widetilde{\Delta}\varphi(z) G(z)dg(z)$ where $G(z)$ is the Green's function for the invariant Laplacian and $dg$ is the volume form associated with the Bergman metric.

We will now prove the estimate
\begin{equation}
\label{2.2}
\int_{\B_n}\abs{f(\la)}^2d\mu(\la)\le \f{(2n)!}{(n!)^2} e \int_{\B_n}\abs{f(z)}^2d\nu(z)
\end{equation}
which will immediately imply the theorem.

First note that
$$
\int_{\B_n}\abs{f(z)}^2d\nu(z)=\f{4n^2n!}{(n+1)\pi^n}\int_{\B_n}\int_{\B_n}\abs{f(z)}^2e^{\varphi(z)}\f{(1-\abs{\lambda}^2)(1-\abs{z}^2)^{n+1}}{\abs{1-\La z,\lambda\Ra}^{2n+2}}G(z)\f{dV(z)}{(1-\abs{z}^2)^{n+1}}d\mu(\lambda).
$$
Using  the estimate $\f{n+1}{4n^2}(1-\abs{z}^2)^n\leq G(z)$ we have 

$$
\int_{\B_n}\abs{f(z)}^2d\nu(z)\geq\f{n!}{\pi^n}\int_{\B_n}\int_{\B_n}\abs{f(z)}^2e^{\varphi(z)}\f{(1-\abs{\lambda}^2)(1-\abs{z}^2)^n}{\abs{1-\La z,\lambda\Ra}^{2n+2}}dV(z)d\mu(\lambda).
$$
\begin{rem} If we did not care about the constant then Theorem \ref{CET-Ball} would be proved.  The reasoning is similar to that in the case of the disc.  
\end{rem}

To obtain the sharper estimate, we again proceed as in the disc.  We focus on the inner integral and will prove the inequality  
\begin{equation}
\label{2.3}
\f{n!}{\pi^n}\int_{\B_n}\abs{f(z)}^2e^{\varphi(z)}\f{(1-\abs{\lambda}^2)(1-\abs{z}^2)^n}{\abs{1-\La z,\lambda\Ra}^{2n+2}}dV(z)\geq \f{(n!)^2}{(2n)!}e^{\varphi(\lambda)}\abs{f(\lambda)}^2\quad\forall\lambda\in\textnormal{supp}\ \mu, 
\end{equation}
which after integration with respect to $d\mu(\la)$ and taking into account that $e^\vf \ge e^{-1}$ gives \eqref{2.2}.  

Consider the conformal change of variables $w=b_\lambda(z)$, where $b_\la$ is an automorphism of the unit ball that exchanges the points $\lambda$ and $0$.  Also observe that $z=b_\la(w)$.  A simple computation shows, see \cite[Theorem 2.2.6]{Rudin} that
$$
dV(w)=\left(\f{1-\abs{\lambda}^2}{\abs{1-\La z,\lambda\Ra}^2}\right)^{n+1}dV(z).
$$
Again following the notation from the previous section, let $\tilde{g}(w):=g\circ b_\lambda(w)$ then the above integral can be recognized as
$$
\f{n!}{\pi^n}\int_{\B_n}\abs{f(z)}^2e^{\varphi(z)}\f{(1-\abs{\lambda}^2)(1-\abs{z}^2)^n}{\abs{1-\La z,\lambda\Ra}^{2n+2}}dV(z)=\f{n!}{\pi^n}\int_{\B_n}e^{\tilde{\varphi}(w)}\abs{\tilde{f}(w)}^2\f{(1-\abs{w}^2)^n}{\abs{1-\La w,\lambda\Ra}^{2n}}dV(w).
$$
In this reduction we have used the algebraic identity  for $b_\lambda$ defined above, namely,
$$
1-\abs{z}^2=\f{(1-\abs{\lambda}^2)(1-\abs{w}^2)}{\abs{1-\La\lambda,w\Ra}^2}.
$$
see \cite[Theorem 2.2.2]{Rudin}.
Continuing the estimate we have 
$$
\f{n!}{\pi^n}\int_{\B_n}e^{\tilde{\varphi}(w)}\abs{\tilde{f}(w)}^2\f{(1-\abs{w}^2)^n}{\abs{1-\La\lambda,w\Ra}^{2n}}dV(w)=\f{n!}{\pi^n}\int_{\B_n}e^{\tilde{\varphi}(w)}\left\vert\f{\tilde{f}(w)}{(1-\La w,\lambda\Ra)^n}\right\vert^2(1-\abs{w}^2)^ndV(w).
$$
The function $\f{\tilde{f}(w)}{(1-\La w,\lambda\Ra)^n}$ is analytic and $\tilde\vf( w)$ is invariant subharmonic (i.e.~$\wt\Delta \wt \vf\ge0$, where $\wt\Delta$ is invariant Laplacian), so (see \eqref{2.1}) the function $u(w) = e^{\tilde{\varphi}(w)}\left| \f{\tilde{f}(w)}{1-\overline{\lambda}w}\right|^2$ is invariant subharmonic.

Integrating in polar coordinates and using the mean value property for invariant subharmonic functions we get 
\begin{eqnarray*}
\f{n!}{\pi^n}\int_{\B_n} u(w) (1-|w|^2)^n d V(w) & = & 2n \int_0^1 (1-r^2)^nr^{2n-1} \int_{\Sp} u(r\theta) d\sigma dr \\
 & \ge & 2n u(0) \int_0^1 (1-r^2)^nr^{2n-1} dr \\
  & = & \frac{(n!)^2}{(2n)!}u(0),
\end{eqnarray*}
where the last integral was recognized as the beta function evaluated at $n+1$ and $n$.
Gathering all together we get 
$$
\f{n!}{\pi^n}\int_{\B_n}e^{\tilde{\varphi}(w)}\abs{\tilde{f}(w)}^2\f{(1-\abs{w}^2)^n}{\abs{1-\La\lambda,w\Ra}^{2n}}dV(w)\geq \f{(n!)^2}{(2n)!} e^{\tilde{\varphi}(0)}\abs{\tilde{f}(0)}^2= \f{(n!)^2}{(2n)!}e^{\vf(\la)}\abs{f(\lambda)}^2.
$$
which is equivalent to \eqref{2.3}. 

This finally shows that 
$$
\int_{\B_n}\abs{f(z)}^2d\mu(z)\leq \f{(2n)!}{(n!)^2}e \norm{f}_{H^2(\B_n)}^2
$$
proving Theorem \ref{CET-Ball} for the ball $\B_n$. \hfill\qed

If one is willing to weaken the initial assumption that 
$$
\sup_{\lambda\in\supp\mu}\int_{\B}\abs{k_{\lambda}(z)}^2d\mu(z),
$$
to instead testing the norm of the reproducing kernels over the support of $\mu$ to testing over all points in the ball $\B_n$, then it is possible to give a slightly different proof of \eqref{2.3}.  One can resort to reproducing kernels for a certain weighted Bergman space to obtain this estimate.


\section{An application to the free interpolation problem}
The classical Carleson Interpolation Theorem says that if the sequence of points $\la_j\in \D$ satisfies the Carleson interpolation condition 
\begin{equation}
\label{CI}
\inf_k \prod_{j\ne k} \left| \frac{\la_k -\la_j}{1- \overline{\la_j} \la_k}\right| =:\delta >0
\tag{C}
\end{equation}
then the sequence $\{\la_j\}_{j=1}^\infty$ is \emph{interpolating}, meaning that for any sequence $\{a_k\}_1^\infty \in \ell^\infty$ there exists a bounded analytic function $f$ such that 
$$
f(\la_k)= a_k, \qquad k=1, 2, \ldots
$$
Moreover, there exists a constant $C$ such that one can always find the interpolating function $f$ satisfying
$$
\|f\|_\infty\le C\|\{a_k\}_1^\infty\|_{\ell^\infty}.
$$
For a long time the only place, where an explicit value of the constant $C=C(\delta)$ was presented, was Nikolskii's book \cite{Nik2}, where it was shown that one can take $C= 32 \delta^{-1} (1+2\ln\delta^{-1})$, see \cite[p.~179]{Nik2}.  A better value of $C$, namely $C= 2e\delta^{-1} (1+2\ln\delta^{-1})$ was given not so long ago by V.~Havin (V.~P.~Khavin) in the appendix of the book \cite{Koosis98} by Koosis. Later in \cite{NikOrSeip2004} the same value of $C$ was obtained by a different method by A.~Nicolau, J.~Ortega-Cerd\`{a} and  K.~Seip. 
Theorem \ref{CET-Disk} gives us another  way to  get the same value  $C= 2e\delta^{-1} (1+2\ln\delta^{-1})$. 

Let us briefly explain how this estimate can be obtained from our result. It was shown in \cite[p.~179]{Nik2} that the constant $C$ can be estimated by $\|\cJ\|\cdot\|\cJ^{-1}\|$, where $J$ is the orthogonalizer of the system $\{k_{\la_j}\}_{j=1}^\infty$. It was also shown there that $\|\cJ\|\cdot\|\cJ^{-1}\|\le \delta^{-1}K^2$, where $K$ is the norm of the embedding operator for the measure $\mu=\sum_k (1-|\la_k|^2) \delta_{\la_k}$. In other words, $K$ is a constant such that for the measure $\mu=\sum_k (1-|\la_k|^2) \delta_{\la_k}$
$$
\int_\D |f|^2d\mu \le K^2 \|f\|, \qquad \forall f\in H^2(\D).   
$$ 

On the other hand, it was also shown in \cite[p.~155]{Nik2} that if the sequence $\{\la_k\}_1^\infty$ satisfies the Carleson condition \eqref{CI}, then for the measure $\mu=\sum_k (1-|\la_k|^2) \delta_{\la_k}$
$$
\sup_{\la\in \D} \int_\D |k_\la|^2 \,d\mu \le 1 + 2\ln\delta^{-1}, 
$$
so by Theorem \ref{CET-Disk} $K^2 \le 2e (1 + 2\ln\delta^{-1})$. The constant in the Carleson Embedding Theorem obtained in \cite{Nik2} was 32, and this accounts for the 32 appearing in Nikolskii's estimate of $C$ for the norm of the operator of interpolation.

\end{document}